\documentclass[a4paper,12pt]{amsart}
\usepackage{amssymb} 
\usepackage{amsmath}
\usepackage{amsthm}
\usepackage{amsbsy}
\usepackage{latexsym}
\usepackage[utf8]{inputenc}
\usepackage[all]{xy}
\usepackage{amsfonts}
\usepackage{mathrsfs}
\usepackage[T1]{fontenc}
\usepackage{mathptmx}
\usepackage{mathbbol}

\usepackage{color}
\definecolor{cobalt}{RGB}{61,89,171}
\usepackage[colorlinks,citecolor=cobalt,linkcolor=cobalt,urlcolor=cobalt,pdfpagemode=UseNone,backref = page
]{hyperref}

\newcommand{\Lie}{\mathscr{L}}      

\theoremstyle{plain}

\newtheorem*{theorem*}{Theorem}

\newtheorem*{conjecture*}{Conjecture}
\theoremstyle{definition}

\theoremstyle{remark}

\newcommand{\bC}{\mathbb{C}}
\newcommand{\bR}{\mathbb{R}}
\newcommand{\R}{\mathbb{R}}

\newcommand{\Z}{\mathbb{Z}}
\newcommand{\fX}{\mathfrak{X}}

\newcommand{\vt}{\vartheta}

\newcommand{\cC}{\mathcal{C}}

\newcommand{\aand}{\quad \textrm{and} \quad}

\title{Locally conformally symplectic and K\"ahler geometry}
\author{Giovanni Bazzoni}

\begin{document}


\begin{abstract}
The goal of this note is to give an introduction to locally conformally symplectic and K\"ahler geometry. In particular, Sections \ref{lcs} and \ref{lck} aim to provide the reader with enough mathematical background to appreciate this kind of geometry. The reference book for locally conformally K\"ahler geometry is \cite{Dragomir-Ornea} by Sorin Dragomir and Liviu Ornea. Many progresses in this field, however, were accomplished after the publication of this book, hence are not contained there -- see the introduction of \cite{Ornea}. On the other hand, there is no book on locally conformally symplectic geometry and many recent advances lie scattered in the literature. Sections \ref{cl_mec} and \ref{susy} would like to demonstrate how these geometries can be used to give precise mathematical formulations to ideas deeply rooted in classical and modern Physics.
\end{abstract}

\maketitle 

\section{Symplectic and locally conformally symplectic geometry}\label{lcs}
A \emph{symplectic manifold} is a smooth manifold $M^{2n}$ with a 2-form $\omega\in\Omega^2(M)$ which is \emph{non-degenerate}, i. e. $\omega^n_p\neq 0$ for every $p\in M$, and \emph{closed}, i. e. $d\omega=0$. The non-degeneracy condition can be rephrased by saying that $\omega$ provides an isomorphism of vector bundles $^\flat\colon TM\to T^*M$, $X\mapsto X^\flat=\imath_X\omega$.

The word \emph{symplectic} was coined by Hermann Weyl in 1939: he replaced the old terminology \emph{complex group} with \emph{symplectic group} to indicate the Lie group of matrices preserving the bilinear skew-symmetric form $\omega_0=\sum_{i=1}^n dx_i\wedge dy_i$ on $\R^{2n}$, see \cite[Page 165]{Weyl}. The etymology is from the Greek $\sigma\upsilon\mu\pi\lambda\varepsilon\kappa\tau\iota\kappa\textrm{\textit{\'o}}\varsigma$, which actually means \emph{complex}.

In the process of getting acquainted with symplectic geometry, something that one experiences quite early is, paraphrasing Gromov, a \emph{curious mixture of ``hard'' and ``soft''}, see \cite{Gromov1} as well as \cite[Page 81]{McDS}. This applies both to the mathematical aspects and to the techniques employed in symplectic geometry. An indication of the soft side of symplectic geometry is certainly \emph{Darboux theorem}, asserting that, locally, two symplectic manifolds can not be distinguished from one another\footnote{This is very different from the Riemannian case, where curvature provides a local invariant.} -- see for instance \cite[Section 8.43]{Arnold} or \cite[Theorem 3.15]{McDS} for a modern proof. Thus, symplectic geometry is somehow a global thing. Of the two conditions ensuring that a 2-form on an even-dimensional manifold is symplectic, however, only one is of global nature, namely closedness. Closedness imposes strong cohomological restrictions on the existence of a symplectic structure on an even-dimensional compact manifold\footnote{The same is not true for open manifolds: as proved by Gromov \cite{Gromov3,Gromov2}, any open manifold with a non-degenerate 2-form admits a symplectic structure.}: for instance, all Betti numbers of even-degree must be non-zero. The general problem of determining which compact manifolds admit a symplectic structure is far from being solved, see \cite{Salamon-D}.

The first true mathematical exposition of what a symplectic manifold is appeared in a paper of Hwa-Chung Lee in 1941, see \cite{Lee}. Lee considers the general setting of an even-dimensional manifold $M^{2n}$ endowed with a non-degenerate 2-form $\omega$. He studies first the \emph{flat} case, in which $d\omega=0$, that is, what is nowadays known as symplectic. Then, he discusses the  problem of two 2-forms $\omega$ and $\omega'$ which are \emph{conformal} to one another: on an open set $U\subset M$ with local coordinates $(x_1,\ldots,x_{2n})$, write
\[
\omega=\sum_{i<j}\omega_{ij}(x)dx_i\wedge dx_j \quad \textrm{and} \quad \omega'=\sum_{i<j}\omega'_{ij}(x)dx_i\wedge dx_j\,;
\]
$\omega$ and $\omega'$ are (locally) conformal if there exists $\varphi\in \cC^\infty(U)$, nowhere vanishing, with $\omega'_{ij}=\varphi \omega_{ij}$. Lee then finds necessary and sufficient conditions for a given $\omega\in\Omega^2(M)$ to be (locally) conformal to a flat, i. e. closed, one: for $n\geq 3$ this happens\footnote{The case $n=1$ is trivial: as remarked by Lee, every $\omega$ is in this case conformal to a flat one, due to dimension reasons. The case $n=2$ is only slightly different -- see the discussion below.} if and only if there exists a 1-form $\vt$ such that $d\omega=\vt\wedge\omega$. It is interesting to notice that the mathematical birthplace of \emph{both} symplectic and locally conformally symplectic geometry is the very same paper of Lee!

The development of symplectic geometry since 1941 was dramatic, kept up first by the French school (Charles Ehresmann, Paulette Libermann, Andr\'e Lichnerowicz, Georges Reeb) in the 1950's, then by the Russian school, with the central figure of Vladimir Arnol'd, and by the American school (Dusa McDuff, Victor Guillemin, Alan Weinstein); another special place is occupied by Mikha\"il Gromov\footnote{The list of quoted mathematicians is of course far from being complete.}. This is however not the right place to extol the ubiquity of symplectic geometry in modern Mathematics -- I refer the reader to the nice surveys \cite{Arnold2,Gotay-Isenberg,McDuff2}.

The fate of locally conformally symplectic geometry, on the contrary, was very different. Except for works of Libermann in 1955 \cite{Libermann1} and Jean Lefebvre in 1966 and 1969 \cite{Lefebvre1,Lefebvre2}, the subject remained in hibernation until the seminal papers of Izu Vaisman: \emph{On locally conformal almost K\"ahler manifolds}, published in 1976 -- see \cite{Vaisman1}, and \emph{Locally conformal symplectic manifolds}, published in 1984 -- see \cite{Vaisman2}.

In \cite{Vaisman1}, Vaisman defines a \emph{locally conformally symplectic manifold}\footnote{Vaisman uses locally \emph{conformal} symplectic, while I stick with the terminology locally \emph{conformally} symplectic in this note. Some recent papers use \emph{conformal symplectic}, see \cite{Chantrain-Murphy,Eliashberg-Murphy}.} as a manifold $M^{2n}$, $n\geq 1$, endowed with a non-degenerate 2-form $\omega\in\Omega^2(M)$ such that every point $p\in M$ has an open neighborhood $U$ such that
\begin{equation}\label{def-lcs-1}
d\left(e^{\sigma}\omega\big|_U\right)=0\,,
\end{equation}
where $\sigma\colon U\to\R$ is a smooth function. If \eqref{def-lcs-1} holds for $U=M$, then $(M,\omega)$ is \emph{globally conformally symplectic}; if it holds for $\sigma$ a constant function, $(M,\omega)$ is clearly a \emph{symplectic} manifold. The work of Lee shows that the above definition is equivalent to the following one: a manifold $M^{2n}$, $n\geq 1$, endowed with a non-degenerate 2-form $\omega\in\Omega^2(M)$, is locally conformally symplectic manifold if there exists a globally defined 1-form $\vt\in\Omega^1(M)$ such that
\begin{equation}\label{def-lcs-2}
d\omega=\vt\wedge\omega \quad \textrm{and} \quad d\vt=0\,.
\end{equation}
The 1-form $\vt$ was baptized the \emph{Lee form} by Vaisman. If $n=1$ one has $d\omega=0=\vt\wedge\omega$ for any choice of $\vt$. For $n\geq 2$, $\vt$ is completely determined by $\omega$; moreover, as remarked by Libermann in \cite{Libermann1}, the second condition in \eqref{def-lcs-2} follows from the first one if $n\geq 3$. $(\omega,\vt)$ is called a locally conformally symplectic structure on $M$. According to this alternative definition, a locally conformally symplectic manifold is globally conformally symplectic if $\vt$ is exact and symplectic if $\vt=0$.

Given a locally conformally symplectic manifold $(M,\omega)$, the \emph{conformal class} of $\omega$ is
\[
\{\omega'\in\Omega^2(M) \mid \exists f\in C^\infty(M) \mid \omega'=e^f\omega\}\,.
\]
If $\vt$ is the Lee form of $(M,\omega)$ and $\omega'=e^f\omega$, then the Lee form of $(M,\omega')$ is $\vt'=\vt+df$, hence the cohomology class of $\vt$ in $H^1_{dR}(M)$ is an invariant of the conformal class.

Formula \eqref{def-lcs-1} implies, in particular, that at a local scale a symplectic manifold can not be distinguished from a locally conformally symplectic manifold. Thus not only all symplectic manifolds locally look alike, in view of Darboux theorem, but potentially there may exist manifolds which locally look like symplectic manifolds and however fail to do so globally! Locally conformally symplectic structures exist on open manifolds, as proved by Fernandes and Frejlich using an $h$-principle -- see \cite{Fernandes-Frejlich}, in particular the Acknowledgements. It was proved very recently by Eliashberg and Murphy using again $h$-principle that a \emph{closed} almost complex manifold $(M,J)$ with a non zero cohomology class $\mu\in H^1(M;\R)$ admits a locally conformally symplectic structure -- see \cite[Theorem 1.8]{Eliashberg-Murphy} for the precise statement. In \cite{Angella-Bazzoni-Parton,Bazzoni-Marrero} explicit examples of compact locally conformally symplectic manifolds which do not admit any symplectic structure are provided.

For this reason, I prefer to consider locally conformally symplectic manifolds as something different from symplectic manifolds. Concretely, this means that our locally conformally symplectic structures will always be assumed to have a Lee form $\vt$ which is not exact.

In his 1976 paper Vaisman proves a few results about locally conformally symplectic manifolds but turns quickly his attention to the metric case, in the wake of Gray's work on almost Hermitian structures. It is in his 1984 article that he extensively studies the non-metric case. Motivated by the metric case, which I will discuss in Section \ref{lck}, Vaisman distinguishes between locally conformally symplectic structures \emph{of the first kind} and \emph{of the second kind}. A locally conformally symplectic structure $(\omega,\vt)$ on $M$ is of the first kind if there exists a vector field $U\in\fX(M)$ such that
\[
\Lie_U\omega=0 \aand \vt(U)=1\,.
\]
Otherwise, it is of the second kind. The above conditions characterize $U$ uniquely; it is the \emph{Lee field} of the locally conformally symplectic structure. A sophisticated way to rephrase this goes as follows: define
\begin{equation}\label{aut:10}
\fX(M,\omega)=\{X\in \fX(M) \mid \Lie_X\omega=0\}\,;
\end{equation}
then $\fX(M,\omega)\subset\fX(M)$ is a subalgebra. If $X\in\fX(M,\omega)$ then $\Lie_X\vt=0$, hence $\vt(X)$ is a constant function on $M$. The \emph{Lee morphism} is $\ell\colon\fX(M,\omega)\to\R$, $\ell(X)=\vt(X)$ and is a morphism of Lie algebras. Thus $(\omega,\vt)$ is of the first kind if and only if the Lee morphism is non zero, hence surjective; of the second kind otherwise. In particular, the Lee form of a locally conformally symplectic structure of the first kind is nowhere zero. I should remark here that in the conformal class of a locally conformally symplectic structure of the first kind there exist always locally conformally symplectic structures of the second kind. To see this, it is enough to choose a function $f$ such that $df_p=-\vt_p$ for some $p\in M$; then the Lee form of $e^f\omega$ will have a zero. Notice that \eqref{aut:10} defines an automorphism of a given element in the conformal class of a locally conformally symplectic structure. If one wants to deal with the whole conformal class, then the object to be considered is the subalgebra
\[
\widehat{\fX}(M,\omega)=\{X\in \fX(M) \mid \exists f_X\in \mathcal{C}^\infty(M) \mid \Lie_X\omega=f_X\omega\}\,;
\]
here $f_X$ should be nowhere 0. In this case as well one sees that the \emph{extended Lee morphism} $\widehat{\ell}\colon\widehat{\fX}(M,\omega)\to\R$, $\widehat{\ell}(X)=\vt(X)+f_X$ is a morphism of Lie algebras (see \cite{Banyaga4}). The Lee morphism and its extended version have been investigated extensively, see for instance \cite{Banyaga4,Haller-Rybicki,Vaisman2}.

Another way to tell locally conformally symplectic structures apart is according to the Morse-Novikov class of the 2-form $\omega$. Given a 1-form $\vt$ on a manifold $M$, one can define a differential operator $d_\vt\colon\Omega^k(M)\to\Omega^{k+1}(M)$ by setting $d_\vt\sigma=d\sigma-\vt\wedge\sigma$. If $\vt$ is closed, then $d_\vt^2=0$ and the \emph{Morse-Novikov}\footnote{The Morse-Novikov cohomology has more than two fathers. In the context of locally conformally symplectic geometry, for instance, it was first considered by Gu\'edira and Lichnerowicz in \cite{Guedira-Lichnerowicz}. It was also considered by Witten in his celebrated paper \cite{Witten}.} cohomology of $(\Omega^*(M),d_\vt)$ is
\[
H^k_\vt(M)=\frac{\ker \{d_\vt\colon \Omega^k(M)\to\Omega^{k+1}(M)\}}{d_\vt(\Omega^{k-1}(M))}\,.
\]

If $M$ is compact, these cohomology spaces are always finite-dimensional, and $H^*_\vt(M)\cong H^*_{dR}(M)$ if $\vt$ is exact. Further, as noticed in \cite{Bande-Kotschick}, the Euler-Poincar\'e characteristic of the Morse-Novikov cohomology of a compact, orientable manifold equals that of the de Rham cohomology, hence it is topological. In general, however, Morse-Novikov cohomology behaves very differently from de Rham cohomology: indeed, if $\vt$ is not exact and $M$ is connected then $H^0_\vt(M)=0$, see \cite{Guedira-Lichnerowicz}; if, in addition, $M^n$ is compact and orientable, then a Poincar\'e duality holds, that is, $H^i_\vt(M)\cong H^{n-i}_\vt(M)^*$, hence $H^n_\vt(M)=0$, see \cite{HallerDiss}.  In \cite{Vaisman3} Vaisman proved that $H^*_\vt(M)$ is isomorphic to the cohomology of $M$ with coefficients in the sheaf of smooth functions $f\in \mathcal{C}^\infty(M)$ which satisfy $d_\vt f=0$. It was proved in \cite{dLLMP} that if $M$ carries a Riemannian metric for which $\vt$ is parallel, then $H^*_\vt(M)=0$. Aside from these general results, the computation of Morse-Novikov cohomology is in general very difficult. For a nilmanifold or a completely solvable solvmanifold\footnote{A \emph{nilmanifold} is the quotient of a connected, simply connected nilpotent Lie group by a discrete and co-compact subgroup. More generally, a \emph{solvmanifold} is a compact quotient of a connected, simply connected solvable Lie group. A solvmanifold is \emph{completely solvable} if the adjoint representation on the corresponding Lie algebra has only real eigenvalues.} the computation of the Morse-Novikov cohomology can be performed algebraically -- see \cite{AndradaOriglia,Angella-Otiman-Tardini,Millionschikov,Mostow}. For more details on the Morse-Novikov cohomology, I refer the reader to \cite{Banyaga,dLLMP,Farber,HallerDiss,Haller-Rybicki}.

The significance of Morse-Novikov cohomology in the context of locally conformally symplectic geometry stems from \eqref{def-lcs-2}: if $(M,\omega,\vt)$ is a locally conformally symplectic manifold then $d\vt=0$ and $d_\vt\omega=d\omega-\vt\wedge\omega=0$, hence the 2-form $\omega$ defines a cohomology class $[\omega]_\vt\in H^2_\vt(M)$. The locally conformally symplectic structure is \emph{exact} if $[\omega]_\vt=0$, \emph{non exact} otherwise. It is easy to see that a locally conformally symplectic structure of the first kind is exact: by defining $\eta=-\imath_U\omega$, where $U$ is the Lee field, one has $\omega=d\eta-\vt\wedge\eta$. The converse does not hold: in fact, being exact is an invariant of the conformal class of a locally conformally symplectic structure, while being of the first kind is not, as it was shown above. The locally conformally symplectic structures constructed by Eliashberg and Murphy in \cite{Eliashberg-Murphy} are exact. The importance of Morse-Novikov cohomology in the context of locally conformally symplectic geometry is highlighted, for instance, by the recent research papers \cite{Angella-Otiman-Tardini,VLV,Otiman}.

Locally conformally symplectic structures of the first kind are strictly related to contact structures. A \emph{(co-orientable) contact structure} on an odd-dimensional manifold $P^{2n+1}$ ($n\geq 1$) consists of a 1-form $\alpha$ such that $\alpha\wedge (d\alpha)^n\neq 0$ at every point, see \cite{Geiges}. Frobenius' integrability theorem shows that the distribution $\xi=\ker\alpha$ is then maximally non-integrable. Let $(P,\alpha)$ be a contact manifold and consider a \emph{strict contactomorphism}, that is, a diffeomorphism $\varphi\colon P\to P$ satisfying $\varphi^*\alpha=\alpha$. Then, as observed for instance by Banyaga in \cite{Banyaga2}, the \emph{mapping torus}\footnote{Given a topological space $X$ and a homeomorphism $\varphi\colon X\to X$, the mapping torus or suspension $X_\varphi$ is the quotient space of $X\times\R$ by the $\Z$-action generated by $1\cdot(x,t)=(\varphi(x),t+1)$. The projection $\pi\colon X_\varphi\to S^1$, $[(x,t)]\mapsto [t]$ is a fiber bundle with fiber $X$. If $M$ is a smooth manifold and $\varphi$ is a diffeomorphism, then $M_\varphi$ is a smooth manifold and $M\to M_\varphi\to S^1$ is a smooth fiber bundle.} $P_\varphi$ admits a locally conformally symplectic structure of the first kind. In the same paper, Banyaga proves a sort of converse to this result: namely, if a compact manifold $M$ is endowed with a locally conformally symplectic structure, then there exist a compact contact manifold $(P,\alpha)$ and a strict contactomorphism $\varphi\colon P\to P$ such that $M$ is diffeomorphic to the mapping torus $P_\varphi$. Banyaga's result, however, does not claim that the original locally conformally symplectic structure on $M$ is the one given by the mapping torus construction. A similar result, in which the given locally conformally symplectic structure is preserved, is proved in \cite{Bazzoni-Marrero2}.

It is interesting to notice that contact and locally conformally symplectic structures come together also in the context of Jacobi structures. According to \cite{Guedira-Lichnerowicz}, indeed, a \emph{transitive} Jacobi manifold is a contact manifold if the dimension is odd and a locally conformally symplectic manifold if it is even.

Locally conformally symplectic structures of the second kind are much less understood.
Concerning, in particular, non exact structures, Banyaga \cite{Banyaga} proved that there exist two families of locally conformally symplectic structures on the 4-dimensional solvmanifold constructed in \cite{ACFM} and that they are non exact. These are the first acknowledged examples of this type of locally conformally symplectic structures. In \cite[Appendix A]{Marrero-Martinez-Padron} it was shown that the locally conformally symplectic structure of the Oeljeklaus-Toma manifolds constructed in \cite{Oeljeklaus-Toma} is not exact. In \cite{Angella-Bazzoni-Parton} the properties of non exact locally conformally symplectic structures extensively are investigated, producing many new examples.

I conclude this section with a collection of results in locally conformally symplectic geometry.

The problem of reduction in locally conformally symplectic geometry was tackled in \cite{Haller-Rybicki,Marrero-Martinez-Padron,Noda}. In \cite{Marrero-Martinez-Padron} the authors also produce \emph{universal models} for exact locally conformally symplectic manifolds, on the line of Tischler's result on universal models for symplectic manifolds, see \cite{Tischler}. A \emph{Moser trick} for locally conformally symplectic forms was proved in \cite{Bande-Kotschick}. The blow-up of a locally conformally symplectic manifold at a point or along a compact symplectic submanifold, i.\, e.\, a submanifold such that the locally conformally symplectic form restricts to a \emph{closed} form, was constructed in \cite{Chen-Yang,Yang-Yang-Zhao}. The notion of \emph{Lagrangian submanifold} makes perfectly sense in the locally conformally symplectic setting. A result on neighbourhoods of Lagrangian submanifolds in locally conformally symplectic manifolds was obtained in \cite{Otiman-Stanciu}, analogous to the known result of Weinstein in the symplectic case, \cite{Weinstein2}. The problem of displacing a Lagrangian submanifold in a locally conformally symplectic manifold is tackled in \cite{Chantrain-Murphy}. The paper also contains some interesting observations on the issues that appear when one tries to apply Floer's machinery or results such as Gromov compactness to the locally conformally symplectic situation. Such issues depend, essentially, on the fact that $\omega$ is not closed, hence no bound \emph{\`a la Gromov} on the energy of a $J$-holomorphic map is possible. The paper \cite{Savelyev} suggests some ideas on how to control the failure of Gromov compactness. The properties of the group of diffeomorphisms preserving the conformal class of a locally conformally symplectic structure are studied in \cite{Haller-Rybicki} -- see also \cite{Banyaga4,Banyaga5,Lefebvre2}. Finally, for a description of locally conformally symplectic structures in the language of Lie algebroids as well as some generalizations I refer the reader to the papers \cite{Kadobianski-Kubarski,KKKW}.

\section{Classical mechanics}\label{cl_mec}

\begin{quote}
Now those Quantities which I consider as gradually and indefinitely increasing, I shall hereafter call \emph{Fluents}, or \emph{Flowing Quantities}, [...] And the Velocities by which every Fluent is increased by its generating Motion, (which I may call \emph{Fluxions}, or simply Velocities or Celerities,) [...]

The Relation of the Flowing Quantities to one another being given, to determine the Relation of their Fluxions.

A relation being proposed, including the Fluxions of Quantities, to find the Relations of those Quantities to one another.
\end{quote}

\begin{flushright}
\emph{Sir Isaac Newton, ``De methodis serierum et fluxionum''}, 1671.\footnote{I am grateful to Prof.\ Antonio Giorgilli for having written amazing lecture notes for the Mathematical Physics courses he taught at University of Milano-Bicocca. As an undergrad I was lucky enough to attend a few of them, a very fruitful experience. His lecture notes contained, among other things, this reference to Newton's original work -- see \texttt{\href{http://www.mat.unimi.it/users/antonio/meccanica/meccanica.html}{http://www.mat.unimi.it/users/antonio/meccanica/meccanica.html}}.}
\end{flushright}

The three sentences of Newton define the objects of interest and summarize the goals of the study of dynamical systems. It was Newton who gave a mathematically precise definition of the three laws that govern classical mechanics, that is, the study of the movement of a body as a response to being exposed to a force. He developed a theory, called in his honor \emph{Newtonian mechanics}, to state and solve the problems posed by classical mechanics, notably arising from planetary motions. In this formalism, the equations of motion of a physical system with $n$ degrees of freedom are given as solutions of $n$ differential equations involving velocities and their derivatives (that is, differential equations of order 2).

Analytical techniques in the study of the problems of classical mechanics, especially celestial mechanics, were brought in by Lagrange at the beginning of the 19\textsuperscript{th} century, founding what is nowadays known as \emph{Lagrangian formalism}; an important role in this formalism is played by the \emph{principle of minimal action}. In particular, as recalled by Weinstein in \cite{Weinstein}, in his 1808 book \emph{M\'emoire sur la th\'eorie des variations des \'el\'ements des plan\`etes}, Lagrange uses explicitly a certain skew-symmetric 6 by 6 matrix. The appearence of geometric techniques in classical mechanics is due to Hamilton, who rewrote Newton's equations as a set of $2n$ differential equations of order 1. In terms of position coordinates $(q_1,\ldots,q_n)$ and corresponding momenta $(p^1,\ldots,p^n)$, the motions are governed by a function $H=H(q_1,\ldots,q_n,p^1,\ldots,p^n)$, the \emph{Hamiltonian} of the system, through the equations
\begin{equation}\label{Hamilton}
\left\{\begin{array}{ccc}
\dot{q}_i & = & \frac{\partial H}{\partial p^i}\\
& &\\
\dot{p}^i & = & -\frac{\partial H}{\partial q_i}\\
\end{array}\right.
\end{equation}
Nowadays\footnote{One could name here many other scientists who contributed to elaborate thorough foundations for classical mechanics -- I prefer to direct the reader to the much more complete references \cite{Abraham-Marsden,Arnold,Cortes-Haupt,Gotay-Isenberg} for further historical and mathematical background.} it is known that, in the simplest case, the phase space of a Hamiltonian system is the cotangent bundle $T^*Q$ of a manifold $Q$ which parametrizes the positions $q$ of the physical system; the corresponding momenta $p$ live on the fibers of the cotangent bundle over a point $q\in Q$ and the Hamiltonian of the system is $H\in C^\infty(T^*Q)$. $T^*Q$ is in a natural way a symplectic manifold; the symplectic form on $T^*Q$ is very easy to describe: if $\pi\colon T^*Q\to Q$ is the canonical projection, define the \emph{Liouville} or \emph{tautological} 1-form $\lambda_{can}\in\Omega^1(T^*Q)$ by $\lambda_{(q,p)}(v)=p(d\pi_{(q,p)}(v))$ for a tangent vector $v$ at $T_{(q,p)}T^*Q$. Then $\omega_{can}=-d\lambda_{can}$ is a symplectic form on $T^*Q$; in local coordinates, one has $\omega_{can}=\sum_{i=1}^ndq_i\wedge dp^i$. In the \emph{Hamiltonian formalism}, the equations of motions \eqref{Hamilton} are given as integral curves of the \emph{Hamiltonian vector field} $X_H$; if $H\colon T^*Q\to \R$ is the Hamiltonian function of the system, then $X_H$ is uniquely determined by the condition $dH=\omega_{can}(X_H,\cdot)$.

From the point of view of Hamiltonian formalism, the fact that non-degeneracy is a local condition implies that the definition of the Hamiltonian vector field is local. Following the illuminating introduction of Vaisman's paper \cite{Vaisman1}, I propose to show that locally conformally symplectic manifolds provide an adequate and more general context for Hamiltonian mechanics. One can make the \emph{Ansatz} that the dynamics on the phase space consists of the orbits of a globally defined vector field $X$. Consider an open set $U_\alpha\subset T^*Q$ with local coordinates $(q_1^\alpha,\ldots,q_n^\alpha,p^1_\alpha,\ldots,p^n_\alpha)$. Then one obtains a local function $H_\alpha\colon U_\alpha\to\R$ such that the orbits of $X$ are defined by a local version of Hamilton's equations,
\begin{equation}\label{Hamilton_local}
\left\{\begin{array}{ccc}
\dot{q}_i^\alpha & = & \frac{\partial H_\alpha}{\partial p^i_\alpha}\\
& &\\
\dot{p}^i_\alpha & = & -\frac{\partial H_\alpha}{\partial q_i^\alpha}\\
\end{array}\right.
\end{equation}
Of course, $X$ is the Hamiltonian vector field of the local Hamiltonian function $H_\alpha$ with respect to the local symplectic form $\omega_{can}^\alpha=\sum_{i=1}^ndq_i^\alpha\wedge dp^i_\alpha$. Suppose $\{U_\alpha\}_\alpha$ is an open covering of $T^*Q$. One the usually requires $\{\omega_{can}^\alpha\}$ and $\{H_\alpha\}$ to piece together to a global symplectic form $\omega_{can}$ and a global Hamiltonian $H$. However, following our \emph{Ansatz}, in order to globalize this local assertion one only needs to prescribe the fact that the transition functions
\begin{equation}\label{canonical}
q_i^\beta=q_i^\beta(q_j^\alpha,p^k_\alpha) \quad \textrm{and} \quad p^i_\beta=p^i_\beta(q_j^\alpha,p^k_\alpha)
\end{equation}
on $U_\alpha\cap U_\beta$ preserve \eqref{Hamilton_local}. Of course, if \eqref{canonical} are \emph{canonical} transformations of the phase space, then $\omega_{can}^\alpha=\omega_{can}^\beta$ and one is back to the symplectic context. However, allowing a \emph{homothetical} change of coordinates, i.\, e.\, taking $H_\beta=\mu_{\beta\alpha}H_\alpha$ for a constant $\mu_{\beta\alpha}\neq 0$, then $\omega_{can}^\alpha=\mu_{\beta\alpha}\omega_{can}^\beta$. Thus our phase space consists of $T^*Q$ with an open covering $\{U_\alpha\}$ and a symplectic form $\omega_{can}^\alpha$ on each $U_\alpha$ such that, on $U_\alpha\cap U_\beta\neq\emptyset$,
\begin{equation}\label{cocycle}
\omega_{can}^\alpha=\mu_{\beta\alpha}\omega_{can}^\beta\,.
\end{equation}
Equation \eqref{cocycle} implies that the collection $\{\mu_{\alpha\beta}\}$ satisfies the \emph{cocycle condition} $\mu_{\beta\gamma}=\mu_{\beta\alpha}\mu_{\alpha\gamma}$, hence one obtains a real line bundle $L\to T^*Q$ with transition functions $\{\mu_{\alpha\beta}\}$. The \emph{global} Hamiltonian is not anymore a smooth function on $T^*Q$ but rather a smooth section of $L$. The cocycle condition can be rephrased by saying that
\[
\mu_{\beta\alpha}=\frac{e^{\sigma_\alpha}}{e^{\sigma_\beta}}
\]
for functions $\sigma_\alpha\colon U_\alpha\to\R$ (resp. $\sigma_\beta\colon U_\beta\to\R$). Now equation \eqref{cocycle} shows that the collection of local 2-forms $\{e^{\sigma_\alpha}\omega_{can}^\alpha\}$ piece together to a global, non-degenerate $2$-form $\omega$ on $T^*Q$. Clearly, the $1$-forms $\{d\sigma_\alpha\}$ piece together to a 1-form $\vt$ and $d\omega=\vt\wedge\omega$. Thus $\vt$ is the Lee form of the locally conformally symplectic structure $(\omega,\vt)$ and $(T^*Q,\omega,\vt)$ is a locally conformally symplectic manifold.

As pointed out in \cite{Haller-Rybicki}, given any manifold $Q$ with a closed 1-form $\bar{\vt}$, the cotangent bundle $T^*Q$ admits a canonical exact locally conformally symplectic structure
\[
(\omega,\vt)=(d_{\vt}(-\lambda_{can}),\vt)\,,
\]
where $\pi\colon T^*Q\to Q$ is the canonical projection and $\vt=\pi^*\bar{\vt}$.\\

I conclude this jaunt into classical mechanics by mentioning a couple of more papers where ideas of conformally symplectic geometry find applications to physical problems.

Let $(M,\omega)$ be a symplectic manifold, let $H\in C^\infty(M)$ be a Hamiltonian function and $X_H$ be the corresponding Hamiltonian vector field. If $f\in C^\infty(M)$ is a function, then the vector field $e^fX_H$ is \emph{conformally Hamiltonian} with Hamiltonian $H$ and conformal factor $e^f$. It clearly satisfies $e^fdH=\omega(e^fX_H,\cdot)$. Moreover, $e^fX_H$ is the Hamiltonian vector field of $H$ for the 2-form $\omega'=e^{-f}\omega$, which is not closed anymore, but conformally closed with Lee form $-df$.

In \cite{MPT} Maciejewski, Przybylska and Tsiganov consider conformally Hamiltonian vector fields in the theory of bi-Hamiltonian systems, in order to produce examples of completely integrable systems.

In \cite{Marle} Marle used conformally Hamiltonian vector fields to study, in a new perspective,
a certain diffeomorphism between the phase space of the Kepler problem and an open subset of the cotangent bundle of $S^3$ (resp. of a 2-sheeted hyperboloid, according to the energy of the motion).

In \cite{Wojtkowski-Liverani}, Wojtkowski and Liverani apply the formalism of conformally symplectic geometry in order to model concrete physical situations such as the Gaussian isokinetic dynamics, also with collisions, and the Nos\'e-Hoover dynamics. More precisely, the authors show that such systems fall under the formalism of conformal Hamiltonian dynamics and explain how to easily deduce results about the symmetric of the Lyapunov spectrum.

\section{K\"ahler and locally conformally K\"ahler geometry}\label{lck}

A K\"ahler manifold is a complex manifold with a compatible Riemannian metric such that the induced complex structure is parallel with respect to the Levi-Civita connection. The Riemannian metric and the complex structure provide a non-degenerate 2-form which is also parallel, in particular closed. Thus K\"ahler geometry lies at the intersection between \emph{complex}, \emph{Riemannian} and \emph{symplectic} geometry. The combination of three geometries produces a class of manifolds which possess distinctive properties within each of the three geometries.

As complex manifolds, K\"ahler manifolds can, to a certain extent, be studied with methods of complex algebraic geometry; indeed, the main source of examples of compact K\"ahler manifolds is provided by projective varieties, i.\, e.\, zero loci of homogeneous polynomials in $\bC P^N$. I should point out, however, that the ``generic'' K\"ahler manifold is not projective. Informally, this question is the content of the \emph{Kodaira problem}:
\begin{quote}
Can every compact K\"ahler manifold be deformed to a projective manifold?
\end{quote}
The answer to this question is, perhaps surprisingly, no, as proved by Claire Voisin in \cite{Voisin1,Voisin2}; see also the survey \cite{Huybrechts2} by Daniel Huybrechts. A certain class of compact K\"ahler manifolds, namely \emph{Hodge manifolds}, can be holomorphically embedded into a complex projective space: this is the content of Kodaira's embedding theorem, see \cite[Theorem 7.11]{Voisin3}. In this case the K\"ahler class is the pullback of the Fubini-Study class but the embedding is, in general, not isometric.

From the perspective of Riemannian geometry, the reduced holonomy of a compact K\"ahler manifold is contained in the unitary group $\mathrm{U}(n)$, where $n$ is half the dimension of the manifold. Manifolds with special holonomy turn out to have many applications in Physics, see for instance \cite{Joyce}.

From the point of view of symplectic geometry, compact K\"ahler manifolds satisfy the \emph{Hard Lefschetz property}, see \cite{Huybrechts1}, while symplectic manifolds need not, see \cite{BFM}. The Lefschetz property implies the well-known fact that the Betti numbers of odd degree are even on a compact K\"ahler manifold (this follows also directly from Hodge theory). In a very actual research area such as homological mirror symmetry, the fact that a symplectic structure is part of a K\"ahler structure on a compact manifold sheds a great deal of light in the study of such duality -- see by way of example \cite{Seidel}. I should also mention here that it was originally believed, and to some extent even erroneously proved, see \cite{Guggenheimer}, that every compact symplectic manifold admitted a K\"ahler metric. It was only in 1976 that Thurston provided the first example of a compact symplectic manifold with first Betti number equal to 3, hence no K\"ahler metric, see \cite{Thurston}. Since then, the quest for compact symplectic manifolds with no K\"ahler metrics has inspired beautiful Mathematics -- see for instance the papers \cite{FM,Gompf,Lupton-Oprea,McDuff} and the book \cite{Oprea-Tralle}.

Finally, concerning the topology of compact K\"ahler manifolds, I should point out that they are formal in the sense of Sullivan, see \cite{DGMS}.

For many purposes\footnote{I will come back to this point in Section \ref{susy}.} it can be convenient to relax the strong integrability properties characterizing the three geometries that come together in a K\"ahler structure. The right framework to do this is that of \emph{almost Hermitian structures}. An almost Hermitian structure on a manifold consists of a triple $(g,J,\omega)$, where $g$ is a Riemannian metric, $J$ is an almost complex structure and $\omega$ is a 2-form, called the K\"ahler form, such that $J$ is an isometry for $g$. Actually two of the three structures determine the third one through the equation
\[
\omega(X,Y)=g(X,JY)\,.
\]
In their celebrated 1980 paper \emph{The sixteen classes of almost Hermitian manifolds and their linear invariants} \cite{Gray-Hervella}, Alfred Gray and Luis Hervella classified almost Hermitian structures in terms of the covariante derivative, with respect to the Levi-Civita connection, of the K\"ahler form. K\"ahler structures are recovered as those almost Hermitian structures whose K\"ahler form is parallel with respect to the Levi-Civita connection. This opens the doors to a whole series of almost Hermitian structures in which some of the integrability properties are not satisfied. Starting with this paper, the study of these structures was undertaken in a systematic way. However, some of them had already appeared before in the literature. For instance, nearly K\"ahler structures were considered by Fukami and Ishihara in 1955 (\cite{Fukami-Ishihara}, on the six sphere) and then studied extensively by Gray \cite{Gray1,Gray2,Gray3}. Locally conformally almost K\"ahler structures were discussed in Vaisman's 1976 paper \cite{Vaisman1}. To an almost Hermitian structure $(g,J,\omega)$ on $M^{2n}$ with $n\geq 2$ one can associate the \emph{Lee form} $\vt\in\Omega^1(M)$, defined as
\[
\vt=-\frac{1}{n-1}J(d^*\omega)\,.
\]
An almost Hermitian structure $(g,J,\omega)$ is \emph{almost K\"ahler} if $d\omega=0$, \emph{locally conformally almost K\"ahler} if $d\omega=\vt\wedge\omega$ and $d\vt=0$. If $J$ is integrable, than the almost K\"ahler structure is K\"ahler and the locally conformally almost K\"ahler is \emph{locally conformally K\"ahler}. For $n=2$, the Lee form of a Hermitian structure is uniquely determined by the condition $d\omega=\vt\wedge\omega$. If $\vt=0$ (resp. $d\vt=0$) then the Hermitian structure is K\"ahler (resp. locally conformally K\"ahkler). If $n\geq 3$ the closedness of the Lee form follows from the equation $d\omega=\vt\wedge\omega$\footnote{This apparent discrepancy between the cases $n=2$ and $n\geq 3$ is due to the fact that, in complex dimension 2, there exist only two ``pure'' classes in the Gray-Hervella classification.}.

Thus, similarly to what happened for K\"ahler manifolds, locally conformally K\"ahler manifolds can be considered simultaneously as complex, Riemannian and locally conformally symplectic manifolds. As I mention above, a K\"ahler manifold $M^{2n}$ can be characterized as a Riemannian manifold whose holonomy lies in $\mathrm{U}(n)$. One could think of a conformal version of manifolds with special holonomy. For instance, locally conformally hyperk\"ahler manifolds are studied in \cite[Chapter 11]{Dragomir-Ornea}. References for locally conformally $\mathrm{G}_2$ and $\mathrm{Spin}(7)$ structures are \cite{FFR,IPP}.

As it happens in the locally conformally symplectic case, a locally conformally K\"ahler manifold is actually K\"ahler, in case $\vt=0$, or \emph{globally} conformal to a K\"ahler manifold if $\vt$ is exact. In general, one can only argue that this conformal property holds locally. I prefer to consider locally conformally K\"ahler manifolds as a class which is distinct from that of K\"ahler manifolds (soemtimes one uses the terminology \emph{strictly} locally conformally K\"ahler). The reference for locally conformally K\"ahler geometry is the monograph \cite{Dragomir-Ornea} by Dragomir and Ornea; see also \cite{Ornea,Ornea-Verbitsky5}.

Of particular importance within locally conformally K\"ahler manifolds are \emph{Vaisman manifolds}\footnote{Vaisman manifolds were first called \emph{generalized Hopf manifolds} by Vaisman, see \cite{Vaisman4}.}; these are characterized by the property that the Lee form is parallel with respect to the Levi-Civita connection. I will implicitly assume that $\|\vt\|\neq 0$ on a compact Vaisman manifold, hence $\vt$ is nowhere zero. An interesting example of a Vaisman \emph{surface} is the Hopf surface. This is defined as a compact complex surface whose universal cover is $\bC^2\setminus \{0\}$. As shown in \cite{Gauduchon-Ornea}, each primary\footnote{A Hopf surface is called \emph{primary} if its fundamental group is isomorphic to $\Z$. Every Hopf surface is finitely covered by a primary one.} Hopf surface admits a locally conformally K\"ahler metric and some Hopf surfaces (those of class 1) admit Vaisman metrics (see also \cite{Parton}). Since every primary Hopf surface is diffeomorphic to $S^3\times S^1$, no Hopf surface admits K\"ahler metrics.

As complex manifolds, locally conformally K\"ahler manifolds are different from K\"ahler manifolds\footnote{In the short note \cite{Aubin}, Aubin erroneously claimed that a compact locally conformally K\"ahler manifold is actually K\"ahler.}. For instance, a small deformation of the complex structure of a K\"ahler manifold remains K\"ahler (see \cite[Theorem 9.23]{Voisin3}). On the other hand, it was shown by Belgun in \cite{Belgun1} that this is not the case for complex structures neither on locally conformally K\"ahler nor on Vaisman manifolds. In the same paper, Belgun carried out a systematical analysis of locally conformally K\"ahler metrics on compact complex surfaces. His paper was in some sense groundbreaking since it was conjectured that all non-K\"ahler compact complex surfaces admitted locally conformally K\"ahler metrics. Locally conformally K\"ahler, albeit non Vaisman, metrics on some Inoue surfaces had been previously constructed by Tricerri in \cite{Tricerri}. Most non-K\"ahler compact complex surfaces admit locally conformally K\"ahler metrics (see \cite{Belgun1,Brunella,Gauduchon-Ornea}); in fact, only one of the three Inoue surfaces is known not to admit any. The \emph{spherical shell conjecture} predicts that every class VII$_0$ surface with $b_2>0$ is a Kato surface, i.\, e.\, it contains a spherical shell\footnote{A \emph{spherical shell} $S$ in a compact complex surface $M$ is a real submanifold diffeomorphic to $S^3$, such that $M\setminus S$ is connected and $S$ has a neighbourhood which is biholomorphic to an annulus in $\bC^2$.}, see \cite{Ornea-Verbitsky8}. If the spherical shell conjecture holds, the remaining non K\"ahler compact complex surfaces admit locally conformally K\"ahler metrics.

A locally conformally K\"ahler manifold can be equivalently defined as a manifold admitting a K\"ahler covering whose deck group acts by conformal transformations (see \cite{Vaisman4}). As proved by Verbitsky in \cite{Verbitsky}, the K\"ahler metric on the universal covering of a Vaisman manifold admits a global K\"ahler potential. Since this property is stable under small deformations, a Vaisman structure deforms to a locally conformally K\"ahler one, not necessarily a Vaisman one. Motivated by this observations, Ornea and Verbitsky defined a class of locally conformally K\"ahler manifolds, which strictly contains Vaisman manifolds, namely \emph{locally conformally K\"ahler manifolds with (proper) potential}, see \cite{Ornea-Verbitsky1,Ornea-Verbitsky2}. Nice results for such manifolds are available. For instance, it was proved in \cite{Ornea-Verbitsky1} that they admit an embedding into a Hopf manifold, provided the complex dimension is at least 3 (see also \cite{Ornea-Verbitsky3}). \emph{Hopf manifolds} are generalizations to arbitrary complex dimensions of Hopf surfaces: they are defined as quotients of $\bC^n\setminus\{0\}$ by a discrete subgroup of linear holomorphisms. A \emph{primary} Hopf manifold is the quotient of $\bC^n\setminus\{0\}$ by the action of the abelian group generated by complex numbers $\lambda_1,\ldots,\lambda_n$, with $0<|\lambda_0|\leq\ldots\leq|\lambda_n|<1$, where the action sends $x_i$ to $\lambda_ix_i$, for $i=1,\ldots,n$ (see \cite{Kamishima-Ornea}). Compact Vaisman manifolds can be embedded into primary Hopf manifolds. In this sense, Vaisman manifolds and, more generally, locally conformally K\"ahler manifolds with proper potential are the analogue of Hodge manifolds in K\"ahler geometry. In locally conformally K\"ahler geometry, the statement corresponding the the Kodaira problem in K\"ahler geometry would be the following:
\begin{quote}
Can every compact locally conformally K\"ahler manifold be deformed to a Vaisman manifold?
\end{quote}

In fact, every compact locally conformally K\"ahler manifold with potential can be deformed to a Vaisman one, as shown in \cite[Theorem 2.1]{Ornea-Verbitsky4}. A compact locally conformally K\"ahler manifold is globally conformally K\"ahler if and only if it admits a K\"ahler metric. Related to the above question, I mention the following two conjectures (see \cite{Vaisman6,Vaisman7}):
\begin{quote}
\begin{itemize}
\item A compact locally conformally K\"ahler manifold satisfying the topological conditions of a K\"ahler manifold admits some global K\"ahler metric.
\item A compact locally but not globally conformally K\"ahler manifold has an odd odd-degree Betti number.
\end{itemize}
\end{quote}
It is easy to see that the first Betti number of a compact Vaisman manifold is odd, hence the second conjecture holds for locally conformally K\"ahler manifolds with potential. A compact complex surface which admits a locally conformally K\"ahler but no K\"ahler metrics has odd first Betti number. In \cite{Oeljeklaus-Toma} Oeljeklaus and Toma disproved the second conjecture by constructing a compact complex 3-fold admitting locally conformally K\"ahler metrics with all odd-degree Betti numbers even. This also settles in the negative the Kodaira problem in the locally conformally K\"ahler context. The so-called \emph{Oeljeklaus-Toma manifolds} are generalizations to arbitrary complex dimensions of Inoue surfaces. They can also be described as solvmanifolds, see \cite{Kasuya}.

The holonomy of locally conformally K\"ahler manifolds has been investigated in \cite{Madani-Moroianu-Pilca}. Although the complex structure is not parallel with respect to the Levi-Civita connection, it can be useful to have an auxiliary metric connection which does fulfill this property. To any Hermitian structure $(g,J)$ on a manifold $M$ one can associate a unique connection, called \emph{Chern connection} $\nabla^C$, which satisfies $\nabla^Cg=0=\nabla^CJ$ and whose torsion $T$ is of type $(2,0)$, that is,
\[
T(JX,Y)=JT(X,Y) \quad \forall X,Y\in\fX(M)\,.
\]
The Chern connection coincides with the Levi-Civita connection if the Hermitian structure is K\"ahler. In \cite{Gauduchon1} Gauduchon associated a 1-form $\tilde{\vt}$ to the Chern connection, the \emph{torsion 1-form}, as follows:
\[
\tilde{\vt}(X)=\textrm{trace}(Y\mapsto T(X,Y))\,.
\]
One can show that $\tilde{\vt}=(n-1)\vt$, hence the Lee form and the torsion 1-form are strictly related. Thus the Lee form of a locally conformally K\"ahler structure measures, in a certain sense, its lack of integrability, where integrability is the K\"ahler case.

A \emph{Weyl structure} on a conformal manifold $(M,c)$ is a torsion-free linear connection $\nabla^W$, the Weyl connection, which preserves the conformal class $c$. This means that there exists a 1-form $\vt$ such that $\nabla^W g=g\otimes\vt$ for every $g\in c$. A conformal Hermitian manifold is a conformal manifold $(M,c)$ with a complex structure $J$ which is Hermitian for some, hence all, $g\in c$. If $\nabla^W J=0$, then $(M,c,J)$ is a \emph{K\"ahler-Weyl manifold}. As pointed out by Kokarev in \cite{Kokarev}, locally conformally K\"ahler manifolds are examples of K\"ahler-Weyl manifolds; the Weyl connection is related to the Levi-Civita connection $\nabla$ by the formula
\[
\nabla^W_XY=\nabla_XY-\frac{1}{2}\vt(X)Y-\frac{1}{2}\vt(Y)X+\frac{1}{2}g(X,Y)U\,,
\]
where $U=\vt^\sharp$ is the \emph{Lee field}. This point of view on locally conformally K\"ahler geometry was adopted in \cite{Kokarev,Kokarev-Kotschick}, with applications to the topology of compact Vaisman\footnote{Kokarev defined in \cite{Kokarev} \emph{pluricanonical} locally conformally K\"ahler metrics (actually K\"ahler-Weyl structures) as those for which $(\nabla\vt)^{1,1}=0$. In \cite{Ornea-Verbitsky4} it was erroneously claimed that a locally conformally K\"ahler metric is pluricanonical if and only if it admits a potential. The mistake was clarified in \cite{Moroianu-Moroianu,Ornea-Verbitsky2}, where it was proved that a compact pluricanonical locally conformally K\"ahler manifold is in fact Vaisman.} manifolds, in particular their fundamental group.

Since the Lee form of a Vaisman structure is parallel, the results of \cite{dLLMP} imply that the underlying locally conformally symplectic structure is exact. But more is true: if $(g,J)$ is Vaisman then, up to a homothety, one can assume that $\|\vt\|=1$ and one can show that the underlying locally conformally symplectic structure is of the first kind, see \cite{Bande-Kotschick,Dragomir-Ornea}; more precisely, one has $\Lie_U\omega=0$ and $\omega=d\eta-\eta\wedge\vt$ for $\eta=-\imath_U\omega$ (see also \cite[Section 9]{Ornea}).

In Section \ref{lcs} we discussed the relation between locally conformally symplectic structures of the first kind and contact structures. A similar relation exists between Vaisman and Sasakian structures. A Sasakian structure is a normal contact metric structure, see \cite{Blair,Boyer-Galicki}. Indeed, the mapping torus of a Sasakian manifold and a Sasakian automorphism, that is, a diffeomorphism which respects the whole Sasakian structure, carries a natural Vaisman structure. In \cite{Ornea-Verbitsky6} the authors claimed that, in the compact case, the converse also holds; as explained in \cite{Ornea-Verbitsky2}, however, the proof is flawed. Nevertheless, the result holds up to diffeomorphism: a compact Vaisman manifold is diffeomorphic to the mapping torus of a Sasakian manifold and a Sasakian automorphism. Morally, this discrepancy between the two directions in similar to what happens in the non-metric case. Based on this approach, a global splitting result for compact Vaisman manifolds was obtained in \cite{Bazzoni-Marrero-Oprea}. As in the non-metric case, let me notice the absence of structure results for compact locally conformally K\"ahler manifolds which are not Vaisman.

Analogous to the symplectic versus K\"ahler case, Ornea and Verbitsky formulated in \cite{Ornea-Verbitsky5} the following problem:
\begin{quote}
Construct a compact locally conformally symplectic manifolds which admits no locally conformally K\"ahler metrics.
\end{quote}
A first answer to this question was provided by Bande and Kotschick in \cite{Bande-Kotschick2}. Different answers are contained in \cite{Bazzoni-Marrero,Bazzoni-Marrero2}.

Related to this problem is a conjecture of Ugarte which aims to give a complete characterization of locally conformally K\"ahler structures on nilmanifolds. In \cite[Page 200]{Ugarte}, he conjectured the following:
\begin{quote}
A compact nilmanifold of dimension $2n\geq 4$ admitting a locally conformally K\"ahler structure is the product of $N$ with $S^1$, where $N$ is a quotient of $H(1,n)$.
\end{quote}
Here $H(1,n)$ is the generalized Heisenberg group, 
\[
H(1,n)=\left\{\begin{pmatrix} 1 & y_1 & y_2 & \ldots & y_n & z\cr
0 & 1 & 0 & \ldots & 0 & x_1\cr
\vdots & 0 & \ddots & \ddots & \vdots & x_2\cr
\vdots & \vdots & \ddots & \ddots & 0 & \vdots\cr
\vdots & \vdots &  & \ddots & 1 & x_n\cr
0 & 0 & \ldots & \ldots & 0 & 1\cr
\end{pmatrix} \mid x_i,y_i,z\in\bR, \ i=1,\ldots,n
\right\}\,.
\]
The conjecture holds in full generality in dimension 4 (\cite{Bazzoni-Marrero}). In higher dimension it holds if one assumes that the complex structure of the locally conformally K\"ahler structure is left-invariant\footnote{This means that it comes from a left-invariant complex structure on the corresponding Lie group.} (see \cite{Sawai}) or if the locally conformally K\"ahler structure is Vaisman (see \cite{Bazzoni}).

We mention here the fact that compact Vaisman manifolds satisfy a Hard Lefschetz property (see \cite{CMDNMY}); this result builds on the Hard Lefschetz property for compact Sasakian manifolds proved in \cite{CMDNY}. Again, the lack of structure theorems for general locally conformally K\"ahler manifolds reflects on the absence of a Hard Lefschetz property in the most general setting.

Compact Vaisman manifolds are, in general, non formal. In 2001, Kotschick introduced the notion of \emph{geometric formality}: a closed manifold is geometrically formal if it admits a Riemannian metric such that the product of two harmonic forms is harmonic (see \cite{Kotschick}). Geometric formality implies formality in the sense of Sullivan, but the converse is not true, see for instance \cite{Kotschick-Terzic}. In \cite{Ornea-Pilca}, Ornea and Pilca showed that geometrically formal compact Vaisman manifolds obey to strong topological restrictions. It is not yet clear the extent to which a compact Vaisman manifold is non formal.

I end this section by quoting some other results about locally conformally K\"ahler manifolds.

Homogeneous locally conformally K\"ahler structures are in fact Vaisman, see \cite{Alekseevsky-Cortes-Hasegawa-Kamishima,Gauduchon-Moroianu-Ornea}. The papers \cite{GOP,GOPP} consider the problem of reduction in locally conformally K\"ahler geometry. In \cite{GOPP} the authors introduce the notions of \emph{presentation} and \emph{rank} of a locally conformally K\"ahler manifold. The rank of a locally conformally K\"ahler structure and his relation with other properties such as the existence of a potential have been further investigated in \cite{Parton-Vuletescu}. Toric locally conformally K\"ahler manifolds, and in particular Vaisman, are considered in \cite{Madani-Moroianu-Pilca2}. The blow-up of a locally conformally K\"ahler manifold was studied in \cite{Ornea-Verbitsky-Vuletescu,Tricerri,Vuletescu}. An interesting contact point between locally conformally symplectic and K\"ahler geometry appears in the papers \cite{Apostolov-Dloussky1,Apostolov-Dloussky2}. The authors consider locally conformally symplectic structures $(\omega,\vt)$ on compact complex surfaces $(M,J)$ such that $\omega$ tames $J$, i.\, e.\, the $(1,1)$-part of $\omega$ is positive definite. The Morse-Novikov cohomology of locally conformally K\"ahler surfaces has been investigated in \cite{Otiman2}. Results on the deformations of Lee classes of locally conformally K\"ahler structures have been obtained in \cite{Goto}. In the more general context of Hermitian structures, metrics which are locally conformal to notable ones, for instance to balanced ones, have been studied in \cite{Angella-Ugarte}.

\section{Strings, Supersymmetry and M-theory}\label{susy}

\begin{quote}
Passer de la m\'ecanique de Newton \`a celle d’Einstein doit \^etre un peu,
pour le math\'ematicien, comme de passer du bon vieux dialecte proven\c{c}al \`a l’argot
parisien dernier cri. Par contre, passer \`a la m\'ecanique quantique, j'imagine, c'est
passer du fran\c{c}ais au chinois.
\end{quote}

\begin{flushright}
\emph{Alexandre Grothendieck, ``R\'ecoltes et Semailles''}, 1986.
\end{flushright}

Given the fact that I am not particularly well-versed in chinese, I will keep this section as low-key as possible. My only goal here is to show that, albeit at a different level from what I discussed in Section \ref{cl_mec}, Physics can motivate and foster research also in the case of locally conformally K\"ahler structures.

At the end of the thirties of the twentieth century, the physical community was able to catch a breath after the establishment of two very important theories, namely General Relativity by Einstein and Quantum Mechanics by Bohr, Heisenberg and Schr\"odinger among others\footnote{Although Einstein is acknowledged as the founding father of General Relativity, his 1905 article ``\"Uber einen die Erzeugung und Verwandlung des Lichtes betreffenden heuristischen Gesichtspunkt'', which won him the Nobel Prize, laid the foundations of the theory of quanta.}. By the end of the seventies the four foundamental forces of Physics, Gravitation, Electromagnetism, Weak and Strong interactions, had been completely described: while General Relativity took care of Gravitation, Quantum Mechanics was able to explain the other three. However, since General Relativity is formulated in the framework of classical physics, in order to fill the gap and elaborate a theory which subsumes the four forces, it is necessary to develop a quantized version of Gravity; this task includes, in particular, the search of a particle, called \emph{graviton}, which carries Gravity. In the Standard Model atoms are broken down to particles, called \emph{fermions} which are the very constituents of matter while the three interactions (electromagnetic, weak and strong) are described as trasmitted by another kind of particles, called \emph{bosons}. These ultimate particles are thought of as punctiform. Among other things, fermions have half-integer spin, while bosons have integer spin. Particle physics in the formalisms of the Standard Model, however, presents some problems.

\emph{String theory} emerged during the sixties and the seventies, with the goal to explain these incongruencies; roughly speaking, in string theory punctiform particles are interpreted as 1-dimensional manifolds with or without boundary. The vibration mode of the string determines the type of particle. In particular, one of the possible states of a string corresponds to the graviton; thus string theory has Quantum Gravity built in. One should keep in mind that, at a large scale, strings look punctiform; the standard analogy with everyday's life is that of a hanging cable or a garden hose: at a certain distance those look one-dimensional, but an ant moving on them would perceive the second dimension.

At a very na\"if level, \emph{Supersymmetry} is a theory in which to each fermion corresponds a boson under a supersymmetry operator. One of the peculiar features of supersymmetry is that it requires a universe with extra dimensions apart from the standard 4-dimensional space time\footnote{The idea of requiring extra dimensions in order to unify Gravity with Electromagnetism goes back to Theodor Kaluza and Oskar Klein in the twenties of last century.}. The number of such extra dimensions is constrained by supersymmetry and can be at most 11 and the universe is thought of as a product $M^n\times K^{11-n}$, where $M$ is Minkowski space-time and $K$ is the so-called \emph{internal space} on which supersymmetry operators act. The extra dimensions encoded in $K$, however, escape our perception. In Physics one uses the term \emph{compactification} to indicate that some extra dimensions ``wrap up'' around a lower dimensional, perceptible universe. This encompasses the idea of string theory that extra dimensions should be ``small''.

Supersymmetry is a further step toward the search for unification and nowadays string theory does include supersymmetry, so that one speaks of \emph{supersymmetric string theory}. Physicists came up with five different supersymmetric string theories: type I, $\textrm{SO}(32)$- and $\textrm{E}_8\times \textrm{E}_8$-heterotic, type IIA and type IIB. Over the years, however, they were able to prove that such theories were related to one another in a highly non-trivial way, through some \emph{dualities}. This led Edward Witten, in 1995, to formulate \emph{M-theory}, a theory which unifies all known supersymmetric string theories; in this sense, all of them are different incarnations of the same theory.

From the point of view of a mathematician, the bridge between M-theory and Geometry is provided by compactification\footnote{See \cite{Ooguri} for a nice explanation of how a string theorist sees geometry.}. Indeed, not every way of compactifying the extra dimensions predicted by M-theory is compatible with the properties of the observed space-time. Indeed, supersymmetry equations for the internal space constrain its geometry. In most models, the internal space is assumed to be a compact Calabi-Yau\footnote{Recall that a Riemannian manifold $(M^{n},g)$ is Calabi-Yau if its holonomy is contained in $\textrm{SU}(\frac{n}{2})$.} manifold. This requirement can be relaxed to include structures with torsion, see for instance \cite{Strominger}.

An explicit compactification of M-theory in 8 dimensions (that is, a solution for which space-time is 3-dimensional and the internal manifold 8-dimensional) was constructed in \cite{Becker-Becker}. The authors show that such 8-dimensional manifold is endowed with a Riemannian metric which, up to a global conformal factor called the warp factor, is Calabi-Yau.

In \cite{Shahbazi}, Shahbazi makes an interesting remark. He asserts that the way in which physicists obtain solutions to their field equations, be they supersymmetry or general relativity, is by performing explicit computations on a local patch of the manifold they are looking for. A (local) solution consists then of an open set with distinguished tensors; the subsequent problem consists in determining which compact manifolds exhibit the particular set as an open set (this is the problem of the maximally analytic extension of a given local patch with a locally defined metric). In order to glue together different open sets on which a solution is known one can require that the change of coordinates respect the distinguished tensors. But more general transformations could be allowed, as we saw in Section \ref{cl_mec}, especially if our goal is to preserve the equations of motions (we refer the reader to \cite{Lazaroiu-Shahbazi} for an explanation of this principle in the setting of supergravity and supersymmetry). In the case we consider here, Shahbazi makes the Ansatz that the warp factor considered in \cite{Becker-Becker} does not necessarily need to be of global nature. If one is able to construct a solution of the equations of motions under this Ansatz, then it would be impossible to distinguish, at a local scale, the two solutions. Shahbazi constructs an explicit compactification of M-theory on a compact 8-dimensional Riemannian manifold whose metric is conformal to a K\"ahler one (actually Calabi-Yau) only \emph{locally}. More precisely, he constructs an internal 8-dimensional locally conformally K\"ahler manifold locally equipped with a preferred Calabi-Yau structure. It is a special type of Hopf manifold, diffeomorphic to $S^1\times S^7$. It is remarkable that the global topology of this explicit compactification does not carry any K\"ahler (hence Calabi-Yau) metric. Shahbazi's solution can be seen as a principal torus bundle over a projective manifold; moreover, the solution is endowed with a codimension 1 foliation whose leaves carry nearly parallel $\textrm{G}_2$-structures. Another nice outcome of Shahbazi's approach is that his solution evades the Maldacena-Nu\~{n}ez \emph{no-go theorem}, see \cite{Maldacena-Nunez}. This asserts, very roughly speaking, that every solution to the equations of supersymmetry which is compatible with a certain zeroth-order approximation of the theory (in a parameter, the Planck length, which corresponds to the tension of the string) must have, in particular, constant warp factor. In order to obtain non-trivial solutions, therefore, one has to allow at least a first-order approximation. In \cite{Becker-Becker} a particular correction of order 6 in the Planck length was included. Shahbazi's solution evades the no-go theorem without needing any kind of correction. The trick lies in the topology of the solution, which is completely different from that of a K\"ahler manifold.

Again we see how a certain relaxation of the K\"ahler (or Calabi-Yau) condition leads to new geometries that can be of use in Physics. It is undeniable that these inputs from Physics are of paramount importance in motivating future research in the area of locally conformally K\"ahler geometry.

\section*{Acknowledgements} 
This note is part of my \emph{Habilitationssarbeit} at the Philipps-Universit\"at Marburg, and most of it was written while I was working there. I am grateful to the whole \emph{Arbeitsgruppe} for the continuous support. The final editing was done in Madrid, where I am supported by a Juan de la Cierva Fellowship at Universidad Complutense de Madrid. I am grateful to Daniele Angella, Marcos Origlia, Alexandra Otiman, Maurizio Parton and S\"onke Rollenske for their comments.

\bibliographystyle{plain}
\bibliography{bibliography}
\end{document}